\documentclass[12pt]{article}
\usepackage[utf8]{inputenc}
\usepackage[english]{babel}
\usepackage{amsfonts}
\usepackage{amssymb}
\usepackage{amsmath}
\usepackage{amsthm}
 \usepackage{CJK}
\usepackage{titlesec}
\titleformat{\section}
  {\normalfont\Large\bfseries\centering}{\thesection}{1em}{}

\numberwithin{equation}{section}
\begin{document}

 \title{Hilbert scheme of rational curves on a generic quintic threefold}
\author{ B. Wang \\
\begin{CJK}{UTF8}{gbsn}
(汪      镔)
\end{CJK}}
\date {Aug, 2018}

\newcommand{\Addresses}{{
  \bigskip
  \footnotesize

   \textsc{Department of Mathematics, Rhode Island college, Providence, 
   RI 02908}\par
  \text{E-mail address}:  \texttt{binwang64319@gmail.com}

}}

\maketitle

\renewcommand{\thefootnote}{\fnsymbol{footnote}} 
\footnotetext{\emph{Key words}: Generic quintic 3-fold, Hilbert scheme, Jacobian matrix,  Clemens'\par\hspace{4.5cc} conjecture} 
   \renewcommand{\thefootnote}{\fnsymbol{footnote}} 
\footnotetext{\emph{2000 Mathematics subject classification }: 14J70, 14J30, 14C05, 14B12}
 
\renewcommand{\thefootnote}{\arabic{footnote}}

\newtheorem{thm}{Theorem}[section]

\newtheorem{proposition}[thm]{\bf {Proposition} }
\newtheorem{theorem}[thm]{\bf {Theorem} }
\newtheorem{ex}[thm]{\bf Example}
\newtheorem{corollary}[thm]{\bf Corollary}
\newtheorem{definition}[thm]{\bf Definition}
\newtheorem{lemma}[thm]{\bf Lemma}

\newcommand{\xdownarrow}[1]{%
  {\left\downarrow\vbox to #1{}\right.\kern-\nulldelimiterspace}
}

\newcommand{\xuparrow}[1]{%
  {\left\uparrow\vbox to #1{}\right.\kern-\nulldelimiterspace}
}

\begin{abstract}  
Let $X_0$ be a generic quintic threefold in projective space 
$\mathbf P^4$ over the complex numbers. For a fixed natural number $d$, 
let $R_d(X_0)$ be  the open sub-scheme of the Hilbert scheme,  
parameterizing irreducible rational curves of degree $d$ on $X_0$.
In this paper, we  show that\par
(1) $R_d(X_0)$ is smooth and of  expected dimension, \par
(2) Combining the Calabi-Yau condition on $X_0$, we further show   \par\hspace{1cc} that it consists of
 immersed rational  curves. \par
 (3) Parts (1) and (2) imply a statement of 
Clemens' conjecture:  \par\hspace{1cc} if  $C_0\in R_d(X_0)$ 
and  $c_0:\mathbf P^1\to C_0$ is the  normalization, the \par\hspace{1cc} normal  sheaf
 is isomorphic to the vector bundle 
$$N_{c_0/X_0}\simeq \mathcal O_{\mathbf P^1}(-1)\oplus \mathcal O_{\mathbf P^1}(-1).$$

\end{abstract}

\tableofcontents

\section{Introduction}

 We work over the complex numbers, $\mathbb C$, and  the Euclidean topology.   A generic (or general)  point is referred  to  a point in a complement of
a proper complex analytic subset of an irreducible analytic (or projective) space.  
The following is the statement of the result.  
\bigskip

\subsection{Statements}

\bigskip

\begin{theorem}  \quad \par
Let $X_0$ be a generic quintic threefold in $\mathbf P^4$ over $\mathbb C$. 
Let $R_d(X_0)$ be  the open sub-scheme of the Hilbert scheme,  
parameterizing   irreducible rational curves of degree $d$ on $X_0$. Then 

(1)  $R_d(X_0)$ is smooth, and of expected dimension
$$ \bigr (dim(X_0)+2-deg(X_0) \bigl)d+dim(X_0)-3. \quad (=0)$$
\par

(2) $R_d(X_0)$  consists of immersed rational curves. Precisely, if $$c_0:\mathbf P^1\to C_0$$ \par\hspace{1cc} is 
the normalization of $C_0\in R_d(X_0)$, then $c_0$ is an immersion, i.e. \par\hspace{1cc} the differential of $c_0$ is injective
everywhere and $c_0$ is an isomorphism \par\hspace{1cc} to its image on an open set.

\end{theorem}

\bigskip

For the same $X_0$, but smooth embedding $c_0$ , Clemens,  in [2], proposed 

\begin{equation} N_{c_0/X_0}\simeq \mathcal O_{\mathbf P^1}(-1)\oplus \mathcal O_{\mathbf P^1}(-1).\end{equation}

\bigskip

Notice that in deformation theory,  part (1) in Theorem 1.1 implies $$H^1(N_{c_0/X_0})=0.$$

Now using part (2) and $X_0$ being Calabi-Yau,  we obtain that the normal sheaf $N_{c_0/X_0}$ must be a vector bundle
decomposed as  
\begin{equation}\begin{array}{cc}
N_{c_0/X_0} \simeq \mathcal O_{\mathbf P^1}(-k-2)\oplus \mathcal O_{\mathbf P^1}(k),\end{array}\end{equation}
with $k\geq -1$.
Applying the Serre duality with the  fact $H^1(N_{c_0/X_0})=0,$
we obtain that $k=-1$. Therefore Theorem 1.1 implies 
\bigskip

\begin{corollary} 
The statement of Clemens' conjecture -- formula (1.1) is correct for all rational curves on $X_0$.
\end{corollary}
\bigskip

\bigskip

{\bf Remark}\quad
  The example of Vainsencher ([5]) shows the immersed, singular rational curves exist on $X_0$. 
 So our corollary proves a statement of the amended Clemens' conjecture ([3]).
\bigskip

\subsection{Outline of the proof}

Our approach  is in the local deformation theory with the focus on algebra, 
 more specifically, 
the linear algebra and the differential algebra.  In  past, the deformation idea  is geometric and is well-known (for instance  see [1]). 
But while it worked in special cases, it  also taught us that in general cases, any kind of deformation of  $X_0$
will encounter the serious hurdle in geometry.  In this paper, we'll  present a method that realizes the  new idea in algebra: \par

 (1) use the down-to-earth approach to set-up the space of rational curves.  \par
(2) then 
  switch the focus  to the first order deformation,  and use linear algebra and differential algebra to reformulate and calculate  
a Jacobian matrix.   

\bigskip

Let's start the problem in geometry.

\bigskip

\subsubsection{ Geometric setting}\par

This is the down-to-earth approach in affine coodinates, that avoids modulo $GL(2)$ group action of $\mathbf P^1$. 
 Let
\begin{equation}M=( H^0(\mathcal O_{\mathbf P^1}(d))^{\oplus 5}\simeq \mathbb C^{5d+5},\end{equation} be
 the space of all
5 tuples of homogeneous polynomials in two variables of degree $d$.
The open set $M_d$ of $M$ represents ( but is not equal to)  $$
\{  c\in Hom_{bir}(\mathbf P^1, \mathbf P^4): deg(c(\mathbf P^1))=d\}.$$
Let $f_0\in H^0(\mathcal O_{\mathbf P^4}(5))$ and $X_0=div(f_0)$ be the 3-fold associated to the quintic section.
  For a generic quintic 3-fold $X_0$, 
    let $I_{X_0}$ be an irreducible component of the incidence scheme
$$ \{  c\in M_d: c^\ast (f_0)=0 \}\subset M_d, $$ 
where we  use the same $c$ to denote  
the regular map $\mathbf P^1\to\mathbf P^4$ associated to the point $c\in M_d$. 
Then we'll prove  
\bigskip

\begin{proposition}\quad\par
$I_{X_0}$ is smooth and  $dim(I_{X_0})$=4.
\end{proposition}

\bigskip

It is well understood that Proposition 1.3 is just the variant  version of the 
main part of Theorem 1.1 ([1], [4]).  But the style of the formulation prepares for a change. 
We switch the focus to another scheme, which will lead us to a different algebra later on.  We call it the projected incidence scheme $I_{\mathbb L}$.
Let 
$$S=\mathbf P( H^0(\mathcal O_{\mathbf P^4}(5))$$ be the space of quintic 3-folds of $\mathbf P^4$.
Let $\mathbb P\subset S$ be a 2-dimensional plane and $\mathbb L\subset  S$ be an analytic open subset.
We define $I_{\mathbb L}$ to be  an irreducible 
component of the scheme
$$
\{  c\in M_d: c^\ast(f)=0,  [f]\in \mathbb L \}\subset M_d, $$
where $[f]$ is a point in $\mathbb L$ representing the quintic $div(f)$. To prove Proposition 1.3, it suffices to prove the 
same property in the local complex geometry 
 for $I_{\mathbb L}$,  
\bigskip

\begin{proposition} \par Assume  $\mathbb L$ satisfies the ``pencil condition"  that  requires all components $I_{\mathbb L}$ 
are disjoint.  Then the following assertions are true.

 \par

(1) For such a sufficiently small subset  $\mathbb L$ of $\mathbb P$,  containing an $S$-generic  \par\hspace{1cc} point,  if all components  
$I_{\mathbb L}$ are smooth and of dimension $6$,  then   \par\hspace{1cc} Proposition 1.3  is true. 
\par (2)  There exist such a plane $\mathbb P$ and its open set $\mathbb L$  containing an $S$-generic  \par \hspace{1cc} $f_0$  such that
  all components $I_{\mathbb L}$ are smooth and of dimension $6$, or  \par \hspace{1cc} equivalently a Jacobian matrix of  each $I_{\mathbb L}$ has full rank.

\end{proposition}

\bigskip

{\bf Remark}\quad 
\par
A more technical description of pencil condition is:
  for the generic $[f_0]\in \mathbb L$, 
 a generic rational curve from  each component $I_{X_0}$ ( with $X_0=div(f_0)$) does not lie in any quintic of $\mathbb L$ other than $X_0$. Because the incidence relation is linear in $S$,  the pencil condition on $\mathbb L$ is also a condition on $\mathbb P$.  The ``pencil condition" is  a necessary condition  for $I_{\mathbb L}$ to be smooth.

\bigskip

\subsubsection{Algebraic setting--the Jacobian matrix}

Theorem 1.1 follows from Proposition 1.3 which follows from   
the part (2) of Proposition 1.4.   In the standard geometry, each $I_{X_0}$ has no difference from
its corresponding $I_{\mathbb L}$.   But we found that  their 
elementary algebra are different.  
Thus we change the focus to algebra and try to prove that

\begin{equation} {\text {    
a Jacobian matrix $J_{\mathbb L}$ of $I_{\mathbb L}$ has full rank. }}\footnote{The same algebraic approach to $I_{\mathbb L}$ 
 fails for 
 $dim(\mathbb L)=0, 1$, and is inaccessible (to us) for $dim(\mathbb L)\geq 3$ .}
\end{equation}

\bigskip

The algebra in (1.4) comes from the affine coordinates of $M$. Let's express $I_{\mathbb L}$ in affine coordinates of $M$. Let $\mathbb P$ be spanned by three non-collinear quintic polynomials $f_0, f_1, f_2$ whose generic open set $\mathbb L$ satisfies the ``pencil condition" and contains an $S$-generic quintic.  Next we define global functions. 
Choose generic $5d+1$ points  $t_i\in \mathbf P^1$ (generic in $Sym^{5d+1}(\mathbf P^1)$), and let
$$\mathbf t=(t_1, \cdots, t_{5d+1})\in Sym^{5d+1}(\mathbf P^1).$$

In the rest of the paper,   we'll also use the following conventions in affine coordinates : \par 
(a) $t_i$ or $t$  denotes a  complex number  in an affine open 
set  $\mathbb C\subset \mathbf P^1$, 
\par (b) for the restriction of $c\in M_d$ to an affine set $\mathbb C\subset \mathbf P^1$, 
 $c(t)$ denotes the \par\hspace{1cc}  image 
$$\begin{array}{ccc} \mathbb C &\stackrel{c}\rightarrow & \mathbb C^5\\
t &\rightarrow & c(t),
\end{array}$$\par

(c) a quintic $f$ is a homogeneous polynomial of degree $5$ in $5$ variables, \par\hspace{1cc}  i.e. $f\in \mathcal O(\mathbb C^5)$.
\par

We should note that in these affine coordinates, the incidence relation $c^\ast(f)=0$ can be expressed
as the composition $f(c(t))=0$ for all $t\in \mathbb C$. 
With these affine expressions, $I_{\mathbb L}$ for an  generic and sufficiently small open set $\mathbb L$ (satisfying the pencil condition) is defined by $5d-1$  polynomials 
\begin{equation} 
\left|  \begin{array}{ccc} 
f_2(c(t_i)) & f_1(c(t_i)) & f_0(c(t_i))\\
f_2(c(t_1)) & f_1(c(t_1)) & f_0(c(t_1))\\
f_2(c(t_2)) & f_1(c(t_2)) & f_0(c(t_2))
\end{array}\right| \end{equation} for $i=3, \cdots, 5d+1$ and variable $c\in M$,  
where $|\cdot |$ denotes the determinant of a 
matrix.  This is because that  the incidence scheme  
$$\{(c, f): f(c(t_i))=0, \ fixed,  \ distinct\ t_1, \cdots, t_{5d+1}\in \mathbb C^1\} $$
implies that $I_{\mathbb L}$ is defined by
\begin{equation} 
\left|  \begin{array}{ccc} 
f_2(c(t_i)) & f_1(c(t_i)) & f_0(c(t_i))\\
f_2(c(t_j)) & f_1(c(t_j)) & f_0(c(t_j))\\
f_2(c(t_k)) & f_1(c(t_k)) & f_0(c(t_k))
\end{array}\right| \end{equation}
for $1\leq i, j, k\leq 5d+1$.
Because of pencil condition, 
$$\begin{array}{c}
\biggl( f_2(c(t_1)), f_1(c(t_1)) ,  f_0(c(t_1))\biggr)\\
\biggl( f_2(c(t_2)),  f_1(c(t_2)),  f_0(c(t_2))\biggr). 
\end{array}$$
are linearly independent for all $c$ in  $I_{\mathbb L}$ for a sufficiently small $\mathbb L$ (by the genericity of $t_1$, $t_2$). This implies 
the polynomials (1.5) and (1.6) define the same scheme around $I_{\mathbb L}$. 
Let $C_M$ be a system of affine coordinates for $M$, which determines an isomorphism 
$$M\simeq \mathbb C^{5d+5}.$$
Then the set of polynomials in (1.5) gives a rise to a holomorphic map
\begin{equation}
M\simeq \mathbb C^{5d+5}\to \mathbb C^{5d-1},\end{equation}
whose differential map is represented by a $(5d-1)\times (5d+5)$ matrix $J_{\mathbb L}$.
 In this paper, we simply call $J_{\mathbb L}$ the Jacobian matrix of the scheme $I_{\mathbb L}$, which  is  uniquely defined 
with the affine coordinates above (We are only concerned with local $J_{\mathbb L}$ around $I_{\mathbb L}$).

\bigskip

\subsubsection{ Jacobian matrices in differential algebra}\par

\bigskip

We introduce differential algebra  to dig deep into the Jacobian matrix $J_{\mathbb L}$.  
We consider differential 1-forms $\phi_i$ on $M$ for $i=3, \cdots, 5d+1$:
  
\begin{equation} \phi_i= \mathbf d
\left|  \begin{array}{ccc} f_2(c(t_i)) & f_1(c(t_i)) & f_0(c(t_i))\\
f_2(c(t_1)) & f_1(c(t_1)) & f_0(c(t_1))\\
f_2(c(t_2)) & f_1(c(t_2)) & f_0(c(t_2))
\end{array}\right|\end{equation}
where $\mathbf d$ is the differential on the regular functions of $M$.  Using these forms, we define a $(5d-1)$-form, 
 \begin{equation} \omega(\mathbb L, \mathbf t)=\wedge_{i=3}^{5d+1}\phi_i \in H^0(\Omega^{5d-1}_M) \end{equation}
 The $\omega(\mathbb L, \mathbf t)$ is the wedge product of 
all the row vectors of  the Jacobian matrix $J_{\mathbb L}$.  
Thus the non-vanishing of it is equivalent to the non-degeneracy of $J_{\mathbb L}$ (of full
rank). 

\bigskip

\begin{proposition} For a generic $\mathbb P$, there exists an open set $\mathbb L$ of the plane $\mathbb P$
consisting of $S$-generic quintics  and satisfying the pencil condition such that 
the form  $\omega(\mathbb L, \mathbf t)$ is nowhere zero on $I_{\mathbb L}$. Therefore $J_{\mathbb L}$ has full rank.

\end{proposition}

\bigskip

 To attack $\omega(\mathbb L, \mathbf t)$, we appeal to differential algebra. We first  expand the determinant (1.5), then 
apply the product rule (from calculus). This yields an expression for $\phi_i$.
\begin{equation}\begin{array}{c}
\phi_i=\left|  \begin{array}{cc} 
 f_1(c(t_1)) & f_0(c(t_1))\\
f_1(c(t_2)) & f_0(c(t_2))
\end{array}\right| \mathbf d f_2(c(t_i))+ \left|  \begin{array}{cc} 
 f_0(c(t_1)) & f_2(c(t_1))\\
f_0(c(t_2)) & f_2(c(t_2))
\end{array}\right|    \mathbf d f_1(c(t_i))
\\+ \left|  \begin{array}{cc} 
 f_2(c(t_1)) & f_1(c(t_1))\\
f_2(c(t_2)) & f_1(c(t_2))
\end{array}\right|\mathbf d f_0(c(t_i))
+ \sum_{l=0, j=1}^{l=2, j=2} h_{lj}^i (c)  \mathbf d f_l(c(t_j)),
\end{array}\end{equation}
where $ h_{lj}^i(c) $ are polynomials in $c$.   We let $c_g$ be a generic point of $\mathbb L$. 
We should notice that each $\phi_i$ at $c_g$ is a linear combination of seven exact 1-forms,
$$
\mathbf d f_3(c(t_i))\big|_{c_g},  \ and \  \mathbf d f_l(c(t_j))\big|_{c_g}, 
$$
where $l=0, 1, 2$, $j=1, 2$,  \begin{equation}f_3=\delta_2 f_2+ \delta_1  f_1+\delta_0 f_0.\end{equation}
for 
$$\delta_2=\left|  \begin{array}{cc} 
 f_1(c_g(t_1)) & f_0(c_g(t_1))\\
f_1(c_g(t_2)) & f_0(c_g(t_2))
\end{array}\right|, \delta_1= \left|  \begin{array}{cc} 
 f_0(c_g(t_1)) & f_2(c_g(t_1))\\
f_0(c_g(t_2)) & f_2(c_g(t_2))
\end{array}\right|, \delta_0=\left|  \begin{array}{cc} 
 f_2(c_g(t_1)) & f_1(c_g(t_1))\\
f_2(c_g(t_2)) & f_1(c_g(t_2))
\end{array}\right|.$$ 
Among these seven  1-forms,  six of them, $\mathbf d f_l(c(t_j))\big|_{c_g}$ are shared by all $\phi_i$.  
Then in differential algebra, the non-vanishing of  $\omega(\mathbb L, \mathbf t)$ evaluated at a point $c_g$ is 
implied by the linear independence of $5d+5=5d-1+6$ differential 1-forms 
\begin{equation}
\mathbf d f_3(c(t_k))\bigg|_{c_g},   \mathbf d f_l(c(t_j))\bigg|_{c_g}, 
\end{equation} in the cotangent space $T^\ast_{c_g}M$, 
for $k=3, \cdots, 5d+1, l=0, 1, 2$, $j=1, 2$.
Notice $5d+5$ is exactly the dimension of $M$. Thus 
we can have the final switch back to a square matrix ( different from the Jacobian matrix $J_{\mathbb L}$), 

\bigskip

\begin{proposition}   The form $\omega(\mathbb L, \mathbf t)$ is non-vanishing at $c_g$ if 
the square  matrix
\begin{equation}\left.
\mathcal A={\partial \biggl(f_3(c(t_3)), \cdots, f_3(c(t_{5d+1})), f_0(c(t_1)), \cdots, f_2(c(t_2))\biggr)\over \partial C_M}\right|_{c_g}
\end{equation}
is non-degenerate, where the row vectors of $\mathcal A$ represent the differential 1-forms in (1.12). 
To emphasize the dependence,  $\mathcal A$ is also expressed as
$$\mathcal A(C_M, f_0, f_1, f_2, \mathbf t).$$

\end{proposition}
\bigskip

\bigskip

\subsubsection{The hidden key to the proof  }

So far we only established the setting, but nothing is substantial.
Next we should see  that the key to the proof hides  inside of  the differential algebra of (1.10). 
 To analyze it, we observe (1.10) carefully.  
First take a pause to note that we are permitted  to fix the generic $f_0$ and manipulate $f_1, f_2$ freely to achieve a non-vanishing $\omega(\mathbb L, \mathbf t)$. 
Among four terms in (1.10), 
the summation $\sum (\cdot)$ and the two other terms with multiples of
differentials 
 \begin{equation}
\mathbf d f_2(c(t_i)),    \mathbf d f_1(c(t_i))\end{equation}
may be accessible by a free choice of $f_1, f_2$.   So the
only  inaccessible term is 
\begin{equation} \left |\begin{array}{cc} 
 f_2(c(t_1)) & f_1(c(t_1))\\
f_2(c(t_2)) & f_1(c(t_2))
\end{array}\right|\mathbf d f_0(c(t_i)). 
\footnote{Proposition 1.3 exactly addressed this lone (without a multiple) differential $d f_0(c(t_i))$ of generic $f_0$ for all $i$.
Actually in local deformation theory,  the differential $d f_0(c(t_i))$ is the only hurdle for Clemens' conjecture. }\end{equation} 

This is due to the genericity of $f_0$. 
But now  we'll let it vanish
by choosing $t_1, t_2$  in
$$ \left |\begin{array}{cc} 
 f_2(c(t_1)) & f_1(c(t_1))\\
f_2(c(t_2)) & f_1(c(t_2))
\end{array}\right|. $$  
So to simplify the hidden algebra in  $\omega(\mathbb L, \mathbf t)$, 
we first fix a generic $c_g$ in $I_{\mathbb L}$, then  arrange
\begin{equation} \left|  \begin{array}{cc} 
 f_2(c_g(t_1)) & f_1(c_g(t_1))\\
f_2(c_g(t_2)) & f_1(c_g(t_2))
\end{array}\right|=0, i.e. \ \delta_0=0.
\end{equation}

In this way we  reduce the $\phi_i$  at $c_g$ to
\begin{equation}
\phi_i|_{c_g}=\biggl ( \mathbf d f_3(c(t_i))+\sum_{l=0, j=1}^{l=2, j=2} h_{lj}^i (c)  \mathbf d f_l(c(t_j))\biggr)\bigg|_{c_g},
\end{equation}
where the troubled term $f_3$ is reduced to $f_3=\delta_1f_1+\delta_2f_2$ which has no constraints. 
The formula (1.17) showed the calculation of $\mathcal A$ boils down  to the calculation of 
differentials $\mathbf d f_3(c(t_i))$ at $c_g$, which can be simplified by specializations. 

\bigskip
 
So key to the proof  in the algebra (1.10)  indicates that it is sufficient to freely  
specialize the first order deformation ( $f_1$ and  $f_2$ )  of the generic quintic $f_0$. This is contrary to that in [1] where  
the whole quintic 3-fold $f_0$ must be deformed to achieve the  same deformation result.

 \bigskip

This key idea,  in terms of the algebra of $\mathcal A$,  is simply the reduction of the matrix. 
It can be achieved by two manipulations: (a) make (1.15) vanish, (b) specialize $f_1, f_2$.
So  to unfold the key idea, we must realize the specializations, to which the most of paper is devoted. 
For this  we choose a specific specialization $f_1, f_2$ (we'll not be surprised to see different kinds of specialization).  
Let
\begin{equation}
f_0=generic, f_1=z_0\cdots z_{2} q, f_2=z_0\cdots z_5
\end{equation}
where $z_0, \cdots, z_4$ are generic homogeneous coordinates of $\mathbf P^4$, 
and $q$ is a  quadratic, homogeneous polynomial in $z_0, \cdots, z_4$.    
Applying a particular type of coordinates, called ``quasi-polar coordinates" $C_M'$ (  section  3.1) associated to 
 the special $f_1, f_2$,  we can break  the Jacobian matrix 
 
\begin{equation}\mathcal A(C_M', f_0, f_1, f_2, \mathbf t')|_{c_g}  \xrightarrow{row\ equiv.}
\left(\begin{array}{cc} D & 0\\
0 & Jac(C_M', c_g)
\end{array}\right)\end{equation}
where $D$ is a non-degenerate diagonal matrix of size $(5d-2)\times (5d-2)$ and $Jac(C_M', c_g)$ is a
$7\times 7$ matrix.  
Next we study the remaining $7\times 7$ matrix,   $Jac(C_M', c_g)$.  
It can be reduced more by specializing $ \mathbf t', c_g, \ and \  q$. 
  Finally we piece them together as a matrix $\mathcal A$,  and deform $$C_M, f_0, f_1, f_2, \mathbf t $$ to  a generic  position. 
  Lastly  $c_g$ can be changed to all rational curves because the number of the rational curves on $X_0$ is finite.

\bigskip

The following arrows give a guide for multiple switches along the logic line (backwards)
\begin{equation}
 R_d(X_0)\Rightarrow I_{X_0}\Rightarrow I_{\mathbb L} \xrightarrow[\text{to algebra}]{\text{from goemtry}}   J_{\mathbb L}
\Rightarrow \omega(\mathbb L, \mathbf t)\Rightarrow \mathcal A
\end{equation}

\bigskip

In the past, the incidence schemes $I_{X_0}$ and $I_{S}$ 
(over the total space $S$) are well studied.   But after $I_{\mathbb L}$ (in (1.20)),  
our ideas of the proof parted ways with the past work in this field.
 
\bigskip
 
\subsection{Technical notations and assumptions}\quad\smallskip

In this subsection, we collect all technical definitions and assumptions used in this paper. Some of them may already be defined before, 
but we'll repeat them more precisely. \bigskip

{\bf Notations}:
\par
(1) $S$ denotes the space all quintics, i.e. $S=\mathbf P(H^0(\mathcal O_{\mathbf P^4}(5)))$.\par
Let $[f]$ denote the image of $f$ under the map
$$\begin{array}{ccc}
H^0(\mathcal O_{\mathbf P^4}(5))\backslash \{0\} &\rightarrow & S. \end{array}$$
\bigskip

(2) Let $$M$$ be
$$(H^0(\mathcal O_{\mathbf P^1}(d))^{\oplus 5}\simeq \mathbb C^{5d+5}$$ 
and $M_{d}$ be the subset that parametrizes  all birational-to-its-image maps $$\mathbf P^1\to \mathbf P^4$$ whose push-forward cycles  have degree $d$. 
\bigskip

(3) Throughout the paper, for $c\in M_d$ the same letter $c$ also denotes its projectivization in $\mathbf P(M_d)$, which is regarded a regular map
 $$c:  \mathbf P^1\to \mathbf P^4. $$
Let  $c^\ast(\sigma)$   denote the pull-back section  of
section $\sigma$ of some bundle over $\mathbf P^4$. The  bundles will not always be specified, but they are apparent in the context.
If $f\in H^0(\mathcal O_{\mathbf P^4}(5))$, $c^\ast(f)$ is also denote by $f(c)$, 
the symbol from the composition in affine coordinates.
If $Y$ is quasi-affine scheme, $\mathcal O(Y)$ denotes the ring of regular functions on $Y$.  

\bigskip

(4) Let $\alpha\in T_{c_0}M$, and $$g: M\to H^0(\mathcal O_{\mathbf P^1}(r))$$
be a regular map. 
Then the  image $g_\ast(\alpha)$ of $\alpha$ under the differential map at $c_0$ is denoted by
the symbol of partial derivatives 
\begin{equation}
{\partial g(c_0(t))\over \partial \alpha}\in T_{g(c_0)}(H^0(\mathcal O_{\mathbf P^1}(r)))=H^0(\mathcal O_{\mathbf P^1}(r)).
\end{equation}
(use the identification $T_{g(c_0)}(H^0(\mathcal O_{\mathbf P^1}(r)))=H^0(\mathcal O_{\mathbf P^1}(r))$).
\bigskip

(5) If $Y$ is a scheme, $|Y|$ denotes the induced reduced scheme of $Y$.

\bigskip

(6)  (a)   If $f\in H^0(\mathcal O_{\mathbf P^4}(5))\backslash \{0\}$ is a quintic polynomial other than $f_0$, we denote the direction of
the line through two points $[f], [f_0]$ in the projective space, $\mathbf P(H^0(\mathcal O_{\mathbf P^4}(5))$  by $\overrightarrow f$.
So $$\overrightarrow f\in T_{[f_0]}\mathbf P(H^0(\mathcal O_{\mathbf P^4}(5)).$$\par
(b) 
Note that the vector $\overrightarrow f$ is well-defined up-to a non-zero multiple. 
In case when $c_0$ can deform to all quintics to the first order, i.e. the map in the formula (2.1) below is surjective, this naturally gives a 
section $<\overrightarrow f>$ of the  bundle $c_0^\ast(T_{\mathbf P^4})$ (may not be unique),  to each deformation $ \overrightarrow f$ of the
quintic $f_0$. This is easily can be understood as the direction of the moving $c_0$ in the deformation $(\overrightarrow f, <\overrightarrow f>)$ of
the pair $(c_0, f_0)$.

\bigskip

(7)  Let $Pr: M_d\times S\to S$ be the projection map. 

Let $\Gamma$  be the union of the open sets of all  irreducible components of the incidence scheme 
\begin{equation}\{(c, [f])\subset M_d\times \mathbf P(H^0(\mathcal O_{\mathbf P^4}(5))): c^\ast(f)=0\}
\end{equation} dominating $S$
such that  fibres  of the map on each component
$$\begin{array}{ccc}\Gamma_i &\rightarrow & Pr(\Gamma_i)\\
\cap & &\cap
\\
\Gamma & & Pr(\Gamma)
\end{array}$$ are equal dimensional. 
  It was shown by Katz ([4]) that there is such a $\Gamma$ with fibre dimension 4. 

\bigskip

(8) Let $f_0, f_1, f_2 \in H^0(\mathcal O_{\mathbf P^4}(5))\backslash \{0\}$ be non-collinear and $f_0$ be generic.
Let $$\mathbb P=span([f_0], [f_1], [f_2])$$  and  
let $\mathbb L$ be an open subset of $\mathbb P$, containing $[f_0]$ and satisfying the pencil condition 

\bigskip

 (9)  Let $\Gamma$ be  that in (7). 

Let
\begin{equation} 
\Gamma_{\mathbb L}\subset \Gamma\cap ( M\times \mathbb L)
\end{equation} be an irreducible component,  and
\begin{equation} 
\Gamma_{f_0}, 
\end{equation} be an irreducible component of
$$P(\Gamma\cap ( M\times \{[f_0]\}))$$
where $P$ is the projection $M\times S \to M$. 
Let  $I_{\mathbb L}=P(\Gamma_{\mathbb L})$ for a component $\Gamma_{\mathbb L}$. In general, these notations  are extended to any subset $B$  of $Pr(\Gamma)\subset S$, so there are 
$I_ B$ and $\Gamma_B$ surjective to $B$.   When $B$ consists of one point $[f_0]$, $I_B$ is also written as 
$I_{f_0}$ or $I_{div(f_0)}$. This notation then is consistent with the notation $I_{X_0}$ used in section 1.2.
Note that $P: \Gamma_{\mathbb L}\to I_{\mathbb L}$ is an isomorphism due to the pencil condition. 

\bigskip

Combing with (7), we obtain  two projections restricted to $\Gamma$, 
$$\begin{array}{ccc}
 & \Gamma &\\
\scriptstyle{Pr}\swarrow &  &\searrow \scriptstyle{P}\\
S & & M.
\end{array}$$
The goal is to know the properties of $Pr$. But we obtain the information through attacking $P$.  
\bigskip

(10)  The term we use the most is ``Jacobian matrix".  
Let $\Delta^n\subset \mathbf C^n$ be an analytic open set with coordinates $x_1, \cdots, x_n$, Let $f_1, \cdots, f_m$ be holomorphic functions on $\Delta^n$. For any positive integers $m'\leq m, n'\leq n$ and a point $p\in \Delta^n$, 
we define
\begin{equation} \begin{array}{c} 
\left.{\partial (f_1, f_2, \cdots, f_{m'})\over \partial (x_1, x_2, \cdots, x_{n'})}\bigr|_{p}: =\left (\begin{array}{cccccc}
{\partial f_1\over \partial x_1}  &{\partial f_1\over \partial x_2}   & \cdots & {\partial f_1\over \partial x_{n'}}   \\
{\partial f_2\over \partial x_1}  &{\partial f_2\over \partial x_2} & \cdots & {\partial f_2\over \partial x_{n'}}  \\
 \vdots & \vdots & \cdots &\vdots \\
{\partial f_{m'}\over \partial x_1}  &{\partial f_{m'}\over \partial x_2} & \cdots & {\partial f_{m'}\over \partial x_{n'}}
\end{array}\right)\right|_p. \end{array}\end{equation}
to be the Jacobian matrix of functions $f_1, \cdots, f_{m'}$ in $x_1, \cdots, x_{n'}$.  So it is a Jacobian matrix with a
particular $C^\infty$ map between Euclidean spaces. But in context we'll skip the descriptions of   Euclidean spaces and
$C^\infty$ maps.
 
\par

This definition coincides with the one used before.

\bigskip
The rest of paper is organized as follows. 

 In section 2, we give a proof of part (1), Proposition 1.4.  It shows some equivalence of 
two different incidence schemes, $I_{X_0}$ and $I_{\mathbb L}$.  
Thus we can shift the focus from $I_{X_0}$ to $I_{\mathbb L}$. These are all in the first order. 
In section 3, we  construct  specializations to calculate differential form 
$\omega(\mathbb L, \mathbf t)$, which represents the Jacobian matrix of $I_{\mathbb L}$.
It leads the proof of Propositions 1.3, 1.4.  Theorem 1.1 is just an invariant 
 expression of them.\bigskip

\bigskip
\section{Equivalence of the incidence schemes, $I_{X_0}$ and $I_{\mathbb L}$}

\subsection{First order deformations of the pair} \quad\smallskip

 Let's start the problem in its first order. It'll lead to some relevant equivalence of $I_{X_0}$ and $I_{\mathbb L}$.

\bigskip

\begin{lemma} 
If $(c_0, [f_0])\in |\Gamma|$ is a generic point, then the projection
\begin{equation} \begin{array}{ccc}T_{(c_0, [f_0])}\Gamma &\rightarrow & T_{[f_0]}S\end{array}\end{equation}
is surjective.

\end{lemma}

\bigskip

\begin{proof} 
Let $|\Gamma |\subset \Gamma$ be the reduced scheme of the scheme $\Gamma$. In a neighborhood of a generic point 
$(c_0, [f_0])\in |\Gamma|$, the projection is a smooth map.
By the assumption,  the projection 
\begin{equation}\begin{array} {ccc} |\Gamma | &\rightarrow &
S\end{array}\end{equation}
is dominant.  Thus 
\begin{equation}\begin{array} {ccc} T_{(c_0, [f_0])}|\Gamma | &\rightarrow &
T_{[f_0]}S\end{array}\end{equation}
is surjective. This proves the lemma

\end{proof}

 To elaborate (6), section 1.3, we apply this lemma to obtain that for any $\alpha\in T_{[f_0]}S$,  there is  
a section denoted by 
$$<\alpha>\in H^0(c_0^\ast(T_{\mathbf P^4}))$$ such that $$ (\alpha, <\alpha>)$$ is tangent to the 
universal hypersurface $$\mathcal X=\{(x, [f])\in \mathbf P^4\times S: x\in div(f)\}.$$ 
 Note that $<\alpha>$ is  unique up to a section in $H^0(c_0^\ast(T_{X_0}))$.  But we will always fix $<\alpha>$ as in introduction.
 
\bigskip

\bigskip
\subsection{The incidence schemes}\quad\smallskip

There are two kinds of incidence schemes $I_{X_0}$ and $I_{\mathbb L}$. In this subsection, we show
certain equivalence between them.  This is meant to shift our focus from $I_{X_0}$ to $I_{\mathbb L}$. 
 More specifically we'll accomplish two goals:
\par
(a) the dimension of the Zariski tangent space of the incidence scheme \par\hspace{1 cc} 
$I_{X_0}$ will force $c_0$ to be an immersion.\par
(b) If $\mathbb L$ satisfies the pencil condition we can reduce the problem to that \par\hspace{1 cc}  over the projected incidence scheme $I_{\mathbb L}$. 

  \bigskip

\begin{lemma} Let $[f_0]\in S$ be a generic point, $\mathbb L_2\subset S$ an open set of the pencil containing 
$f_0$ and another  quintic $f_2$.   Assume they determine the components $ I_{f_0},  I_{\mathbb L_2}$ satisfying
\begin{equation}
I_{f_0}\subset I_{\mathbb L_2}, c_0^\ast (f_2)\neq 0 \ for\ generic\ \ c_0\in I_{f_0}.
\end{equation}
 Then
 
\par
(a)  \begin{equation}
{T_{c_0}I_{f_0}\over ker} \simeq H^0(c_0^\ast(T_{X_0})).\end{equation}
where $ker$ is a line in $T_{c_0}I_{f_0}$.  

\par
(b) 
\begin{equation} dim( T_{(c_0, [f_0])}\Gamma_{\mathbb L_2})=dim( T_{c_0}I_{f_0})+1, 
\end{equation}
and furthermore
\begin{equation} dim( T_{c_0}I_{\mathbb L_2}))=dim( T_{c_0}I_{f_0})+1, 
\end{equation}
\par

(c) If $dim( T_{c_0}P(\Gamma_{\mathbb L_2}))$=5, 
then\par

(1) $c_0$ is an immersion, \par
 (2) and 
\begin{equation}
H^1(N_{c_0/X_0})=0.\end{equation}
\end{lemma}

\bigskip

\begin{proof}
(a). Let $a_i(c, f), i=0, \cdots, 5d$ be the coefficients of 
polynomial $f(c(t))$ in parameter $t$. Then the scheme
$$\Gamma$$ in $M\times \mathbf P^4$ is defined by homogeneous
polynomials $$a_i(c, f)=0, i=0, \cdots, 5d, \ locally\ around\ (c_0, [f_0]).$$ 
 Let $\alpha\in T_{c_0}M$. The equations on $\alpha$
\begin{equation} {\partial  a_i(c_0, f_0)\over \partial \alpha}=0, \ all\ i \end{equation}  by the definition,  are necessary and sufficient conditions
for $\alpha$ to lie in $$T_{(c_0, [f_0])}I_{f_0}.$$ 
On the other hand there is an evaluation map $e$:
\begin{equation}\begin{array}{ccc}
M\times \mathbf P^1 &\rightarrow & \mathbf P^4\\
(c, t) &\rightarrow & c(t)
\end{array}\end{equation}
The differential map
$e_\ast$ gives a morphism $e_m$:
\begin{equation}\begin{array}{ccc}
T_{c_0}M &\stackrel{e_m}\rightarrow & H^0(c_0^\ast( T_{\mathbf P^4}))\\
\alpha &\rightarrow & e_\ast (\alpha)
\end{array}\end{equation}
(Note $c_0$ is birational to its image. Thus $c_0^\ast( T_{\mathbf P^4})$ exists).
Suppose there is an $\alpha$ such that
$e_\ast (\alpha)=0$. We may assume $c_0$ is a map
$$\mathbb C^1\to \mathbb C^5\backslash\{0\}.$$
 
Since $c_0$ is birational to its image, there is a Zariski open set  $$U_{\mathbf P^1}\subset \mathbb C^1\subset \mathbf P^1$$
and an open set $$V\subset c_0(\mathbf P^1)\subset \mathbb C^{5}\backslash\{0\}$$
such that $c_0|_{U_{\mathbf P^1}}$ is an isomorphism
$$U_{\mathbf P^1} \to V.$$ 
Due to the equation $e_\ast (\alpha)=0$, on $T_{t}U_{\mathbf P^1}$
$$(\alpha_0(t), \cdots, \alpha_4(t))=\lambda(t) c_0(t)$$
on $V$ (at each point $(c_0(t), \cdots, c_4(t))$ of $V$)
where $\lambda(t)$ lies in $\mathcal O(U_{\mathbf P^1})$. Because
$(\alpha_0(t), \cdots, \alpha_4(t))$ is parallel to
$c_0(t)$ at all points $t\in \mathbf P^1$, $\lambda(t)$ can be extended to $\mathbf P^1$. Hence
$\lambda(t)$ is in $H^0(\mathcal O_{\mathbf P^1})$. So it is a constant (independent of $t$).
Therefore $\alpha\in \mathbb C^{5d+5}$ is parallel to $$c_0\neq 0\in \mathbb C^{5d+5}.$$
This shows that 
$$dim(ker(e_m))=1.$$
  By the dimension count, $e_m$ must be surjective.
 
For any $\alpha\in H^0(c_0^\ast( T_{\mathbf P^4}))$, 
$\alpha\in H^0(c_0^\ast(T_{X_0}))$ if and only if
\begin{equation} {\partial  f_0(c_0(t))\over \partial \alpha}=0, \end{equation} 
for generic $t\in \mathbf P^1$. Notice equations (2.9) and (2.12) are exactly the same. Therefore
$e_m$ induces an isomorphism 
\begin{equation} \begin{array}{ccc} {T_{c_0}I_{f_0}\over ker(e_m)} &\stackrel{e_m}\rightarrow &  H^0(c_0^\ast(T_{X_0}))
\end{array} \end{equation}
This proves part (a).\par

(b).  Let \begin{equation}
<\overrightarrow {f_2}>_M\in e_m^{-1} (<\overrightarrow {f_2}>)\end{equation}
 be an inverse of the
the section $<\overrightarrow {f_2}>$ in the map (2.11). 
Since \begin{equation}
{\partial f_0(c_0(t))\over \partial <\overrightarrow {f_2}>}=-c_0^\ast(f_2)\neq 0,
\end{equation}
by part (a),  $<\overrightarrow {f_2}>_M$ does not lie in $T_{c_0}I_{f_0}$. Hence
\begin{equation}
T_{c_0}P(\Gamma_{\mathbb L_2})=T_{c_0}I_{f_0}+ \mathbb C<\overrightarrow {f_2}>_M
\end{equation}
has dimension $dim(T_{c_0}I_{f_0})+1$, where $P$ is the projection 
 $$\begin{array} {ccc} P: M\times S &\rightarrow & M \end{array}$$

Notice $$T_{c_0}P(\Gamma_{\mathbb L_2})\simeq T_{(c_0, f_0)}\Gamma_{\mathbb L_2}.$$
Hence \begin{equation}
dim(T_{(c_0, f_0)}\Gamma_{\mathbb L_2})=dim(T_{c_0}I_{\mathbb L_2})=dim(T_{c_0}I_{f_0})+1.
\end{equation}
\bigskip

(c) If  $dim( T_{c_0}P(\Gamma_{\mathbb L_2}))$=5, then by part (a), (b), 
\begin{equation} 
dim(H^0(c_0^\ast(T_{X_0})))=3.
\end{equation}
Now we consider it from a different point of view.
Because $c_0$ is a birational map to its image,  there are finitely many points $t_i\in\mathbf P^1$ where the differential map

\begin{equation}\begin{array}{ccc}
(c_0)_\ast: T_{t_i}\mathbf P^1 &\rightarrow & T_{c_0(t_i)}\mathbf P^4
\end{array}\end{equation}
is a zero map. Assume its vanishing order at $t_i$ is $m_i$ .  Let
\begin{equation}
m=\sum_i m_i.
\end{equation}
Let $s(t)\in H^0(\mathcal O_{\mathbf P^1}(m))$ such that $$div(s(t))=\Sigma_i m_it_i.$$

The sheaf morphism $(c_0)_\ast$ is injective and induces a composed morphism $\xi_s$ of sheaves
\begin{equation}\begin{array}{ccccc}
T_{\mathbf P^1} &\stackrel{(c_0)_\ast}  \rightarrow & c_0^\ast(T_{X_0}) & \stackrel{1\over s(t)}\rightarrow &
c_0^\ast(T_{X_0})\otimes \mathcal O_{\mathbf P^1}(-m).\end{array}
\end{equation}

It is easy to see that the induced bundle morphism $\xi_s$ is injective. Let 
\begin{equation}
N_m={c_0^\ast(T_{X_0})\otimes \mathcal O_{\mathbf P^1}(-m)\over \xi_s(T_{\mathbf P^1})}.
\end{equation}
Then
\begin{equation}
dim(H^0(N_m))=dim(H^0(c_0^\ast(T_{X_0})\otimes \mathcal O_{\mathbf P^1}(-m)))-3.
\end{equation}
On the other hand, three dimensional automorphism group of $\mathbf P^1$ gives a rise
to a 3-dimensional subspace  $Au$ of  $$H^0(c_0^\ast(T_{X_0})).$$
By (2.18), $Au=H^0(c_0^\ast(T_{X_0}))$. Over each point $t\in \mathbf P^1$, $Au$ spans 
a one dimensional subspace. Hence 
\begin{equation} c_0^\ast(T_{X_0})=\mathcal O_{\mathbf P^1}(2)\oplus \mathcal O_{\mathbf P^1}(-k_1)\oplus \mathcal O_{\mathbf P^1}(-k_2),
\end{equation}
where $k_1, k_2$ are some positive integers.
This implies that
\begin{equation} dim (H^0(c_0^\ast(T_{X_0})\otimes \mathcal O_{\mathbf P^1}(-m)))=dim( H^0(\mathcal O_{\mathbf P^1}(2-m)).
\end{equation}
Then
\begin{equation} dim (H^0(c_0^\ast(T_{X_0})\otimes \mathcal O_{\mathbf P^1}(-m)))=3-m.\end{equation}

Since $dim(H^0(N_m))\geq 0$, by the formula (2.23), $-m\geq 0$. 
By the definition of $m$, $m=0$.  Hence $c_0$ is an immersion. 
\par
Next we prove (2). Notice that $(c_0)_\ast (T_{\mathbf P^1})$ is a sub-bundle generated by
global sections. It must be the $\mathbf O_{\mathbf P^1}(2)$ summand in (2.24) because $k_1, k_2$ are positive.
Therefore 
\begin{equation}
N_{c_0/X_0}\simeq \mathcal O_{\mathbf P^1}(-k_1)\oplus \mathcal O_{\mathbf P^1}(-k_2).
\end{equation}
Since $deg(c_0^\ast(T_{X_0}))=0$,  $k_1=k_2=1$.

Therefore 
\begin{equation}
H^1(N_{c_0/X_0})=0.\end{equation}

\end{proof}

\bigskip

Now we can describe the case for $\mathbb P$. 
Recall $\mathbb L$ is an open set of $\mathbb P$ spanned  by $f_0, f_1, f_2$, where $f_0$ is generic. 
 Let  $$\mathbb L_2\subset \mathbb L$$ be a pencil 
containing the generic $[f_0]$. 
\bigskip

\begin{lemma} Let $[f_0], [f_1], [f_2]$ be non-collinear quintic hypersurfaces and $f_0$ is generic. Let $\mathbb L_2$ be an open set of
the pencil $span([f_0], [f_2])$ as in Lemma 2.2. 
Also assume that  an open set $\mathbb L$ of the span of $[f_0], [f_1], [f_2]$ satisfies the pencil condition.
 and $$f_0\in \mathbb L_2\subset \mathbb L$$ as before.  We choose components
$$I_{f_0}\subset I_{\mathbb L_2}\subset I_{\mathbb L},$$
and let $c_0\in I_{f_0}$. 
Then
\begin{equation} 
dim( T_{c_0} I_{\mathbb L})=dim( T_{c_0}I_{\mathbb L_2})+1. 
\end{equation}

\end{lemma}

\bigskip

\begin{proof} The formula (2.16) asserts  \begin{equation}
 T_{c_0}I_{\mathbb L_2}\simeq T_{c_0}I_{f_0}\oplus \mathbb C<\overrightarrow{f_2}>_M.\end{equation}

Note \begin{equation}
T_{c_0}I_{\mathbb L}\simeq T_{c_0}I_{f_0}\oplus \mathbb C<\overrightarrow{f_2}>_M+ \mathbb C<\overrightarrow{f_1}>_M.
\end{equation}

Thus it suffices to show 
the section $<\overrightarrow{f_1}>_M$ does not lie in $$T_{c_0}I_{f_0}\oplus \mathbb C<\overrightarrow{f_2}>_M.$$ 
By the  notation (6) in section 1.3, this is equivalent to 
$$ f_1(c_0(t))+\epsilon_2 f_2(c_0(t))\neq 0$$ for all complex numbers $\epsilon_2$. By the 
``pencil condition", the  lemma is proved.

\end{proof}

\bigskip

Lemmas 2.2, 2.3 imply\bigskip

\begin{corollary}
 The part (1), Proposition 1.4 is correct.

\end{corollary}
\bigskip

So in the rest paper, we concentrate on part (2) of Proposition 1.4.\bigskip
\bigskip

\section{Projected incidence scheme $I_{\mathbb L}$}
In this section, all neighborhoods and the word ``local" are in the sense of Euclidean topology. 
The  topic of this section is the computation of the differential from
$$\omega(\mathbb L, \mathbf t)$$
through multiple block Jacobian matrices of $\mathcal A$ under specializations. This mainly corresponds to
the ``specialization" step in the key to the proof in section 1.2.4. 
Specifically  there are only two methods in the calculations of matrices (they must be specialized first):\par
(1) For a Jacobian matrix of large size, we construct special coordinates  \par\hspace{1cc} -- quasi-polar coordinates, 
to reduces it to a diagonal matrix.
\par
(2) For a Jacobian matrix of smaller size, use the free choice of $f_1, f_2$ and   \par\hspace{1cc} other parameters 
 to  show the linear  independence of
row vectors. \bigskip

 We organize them in three subsections

\bigskip

(3.1) Construct local analytic coordinates of $M_d$, that are associated  \par\hspace{2cc} to the specialization of the quintics $f_1, f_2$.\par
(3.2) Use above coordinates to show the largest block matrix of $\mathcal A$  \par\hspace{2cc}  (in (1.19)),   associated to  $\omega(\mathbb L, \mathbf t)$  is a diagonal  matrix. \par
(3.3) Put all preparations above in the context to construct and   compute \par\hspace{2cc}$ \omega(\mathbb L, \mathbf t)$.

\bigskip

\subsection{Quasi-polar coordinates}

\bigskip

We introduce local analytic coordinates of the affine space $M$, that will simplify expressions 
of differentials on $M$. \bigskip 

\begin{definition} (polar coordinates)
Let $a_0\in H^0(\mathcal O_{\mathbf P^1}(d))$ be a non-zero element satisfying the zeros are distinct.
Then there is a Euclidean neighborhood $U\subset H^0(\mathcal O_{\mathbf P^1}(d))$ of $a_0$, which
has analytic coordinates $r, w_1, \cdots, w_d$ ($r\neq 0$) such that
all elements $a\in U$ has an expression
\begin{equation}
a=r \prod _{j=1}^d  (t-w_j).
\end{equation}
We call $\{r, w_1, \cdots, w_d\}$ the polar coordinates of  $H^0(\mathcal O_{\mathbf P^1}(d))$ at $a_0$.
\end{definition}
\bigskip

Next we  fix the notations of polar coordinates in the $M_d$. 
 Let non-zero $ c=( c^0, \cdots,  c^4)$ 
with $$ c^i\in H^0(\mathcal O_{\mathbf P^1}(d)), i=0, \cdots, 4$$
be a varied point of $M_d$ in a small analytic neighborhood. 
We assume the equations $ c^i(t)=0, i\leq 4, $ always have  $5d$ distinct  zeros $$ \theta_j^i, for\ i\leq 4, j\leq d$$
The we have polar coordinates for $M$ (around some points). 
We denote them by
\begin{equation}\begin{array}{c}
r_i, \theta_j^i,\\
 j=1, \cdots, d, i=0, \cdots, 4\end{array}\end{equation}
with $ r_i\neq 0$  satisfying
\begin{equation} c^i(t)=r_i \prod _{j=1}^d (t-\theta_j^i).  
\end{equation}

The values of the center point $c_g$ of the neighborhood  are denoted by
$$\begin{array}{c} \mathring{r_i}, \mathring{ \theta_i^j} \\
\ for\ i=0, \cdots, 4, j=1, \cdots, d.
\end{array}$$
\bigskip

Next we define quasi-polar coordinates that are associated to the special quintics we are going to
choose later.  They are the polar coordinates for $M$ with a deformation for  last two components $c^3, c^4$.
 Let $q$ be a homogeneous quadratic polynomial in variables $z_0, \cdots, z_4$.
Let 
\begin{equation} 
h(c, t)=\delta_1q(c(t))+\delta_2 c^3(t) c^4(t).
\end{equation}
for $c\in M$, where $\delta_i, i=1, 2$ are two  complex numbers, generic in $\mathbb C^2$. 
Assume for $c$ in a small analytic neighborhood, $h(c, t)=0$ has $2d$ distinct zeros.
   Let $\gamma_1, \cdots, \gamma_{2d}$ be the zeros of 
$h(c, t)=0$. Similar to the polar coordinates, we let  
$$h(c, t)=R \prod _{k=1}^{2d}(t-\gamma_k), R\neq 0 $$
It is clear that $$\begin{array}{c}
R=\delta_1 q(r_0, r_1, r_2, r_3, r_4)+\delta_2 r_3r_4, \ and\\
\gamma_k\ are \ analytic\ functions\ of\ c.
\end{array}$$
( Notice $R$ is the value of $h(c, t)$ at $t=\infty$, the coefficient of the highest order.).
Let the coordinates values at the center point be $\mathring{R}, \mathring{\gamma_k}$.

\bigskip

\begin{proposition} Let $(\delta_1, \delta_2)$ and $q$ be generic. 
Let $U_{c_g}\subset M$ be an analytic neighborhood of a center point $c_g$ as above. 

\par
 Let \begin{equation}\begin{array}{ccc}
\varrho: U_{c_g} &\rightarrow &  \mathbb C^{5d+5}
\end{array}\end{equation}
be a regular map that is defined
by 
\begin{equation}\begin{array}{cc} &
(\theta_0^1, \cdots, \theta_4^d, r_0, r_1, r_2, r_3, r_4) \\
 &\xdownarrow{0.3cm}\scriptstyle{\varrho}\\ &
(\theta_0^1, \cdots, \theta_2^d, r_0, r_1, r_2, r_3, r_4,  \gamma_1, \cdots, \gamma_{2d}).
\end{array}\end{equation} 
Then $\varrho$ is an isomorphism to its image.
\end{proposition}
\bigskip

\begin{proof}
It suffices to prove the complex differential of $\varrho$ at $c_g$ is an isomorphism for a specific $q$, $\delta_i$.  So we assume  that 
$$\delta_1=0, \delta_2=1.$$
Then $h(c, t)=c^3(t)c^4(t)$. Hence $\gamma_k, k=1, \cdots, 2d$ are just $$\theta_j^i, i=3, 4, j=1, \cdots, d.$$
 So $\varrho$ is the identity map.
We complete the proof. 

\end{proof}

\bigskip

\begin{definition}
By Proposition 3.2,  for generic $(\delta_1, \delta_2), q$, 
\begin{equation} \theta_0^1, \cdots, \theta_{2}^d, r_0, r_1, r_2, r_3, r_4, \gamma_1, \cdots, \gamma_{2d}
\end{equation} 
are  local analytic coordinates of $M$ around $c_g$, and
$c_g$ corresponds to the coordinate values with $\circ$ accent.  We denote the system of coordinates by 
$$C_M'$$
and will be called quasi-polar coordinates. 
\bigskip

\end{definition}
\bigskip

\subsection{The largest block matrix} 

The specifically tailored  quasi-polar coordinates above will automatically imply 
that  the largest block in (1.19) is   a diagonal matrix.  \bigskip

 Let's define this matrix. Choose  a generic homogeneous coordinates  $[z_0, \cdots, z_4]$ for $\mathbf P^4$. 
Let 
\begin{equation} 
f_3=z_0z_1z_2(\delta_1q+ \delta_2z_3 z_4).
\end{equation}
be a quintic polynomial, 
where $(\delta_1, \delta_2), q$ are generic.
(Later we will choose $f_1=z_0\cdots z_4, f_2=z_0z_1z_2q$ for the specialization).
Let 
$c_g\in M_d$ such that 
$$f_3(c_g(t))\neq 0$$
Recall we have  denoted the zeros of $c_g^i(t)=0$ by $\mathring{\theta_j^i}$ and zeros of
\begin{equation}
(\delta_1q+ \delta_2z_3 z_4)|_{c_g(t)}=0
\end{equation} for varied $c$
by $\gamma_k, k=1, \cdots, 2d$. 
We assume $\mathring{ \theta_j^i}, i=0, \cdots, 4, j=1, \cdots, d$ are distinct, and $\mathring{\gamma_k}, k=1, \cdots, 2d$ are also distinct. 
For the simplicity we denote $5d$ complex numbers
$$\mathring{\theta_1^0},\cdots, \mathring{\theta_2^d}, \mathring{\gamma_1}, \cdots, \mathring{\gamma_{2d}}$$
by 
$$\tilde t_1, \tilde t_2, \cdots, \tilde t_{5d}.$$

So $\tilde t_1, \tilde t_2, \cdots, \tilde t_{5d}$ are zeros of $f_3(c_g(t))=0$.

\bigskip

\begin{lemma}   Recall 
in Definition 3.3, 
$$ \theta_0^1, \cdots, \theta_2^d, \gamma_1, \cdots, \gamma_{2d}, r_0, \cdots, r_4  
$$ are analytic coordinates of $M$ around the point $c_g$.

Then \par

(a) the Jacobian matrix

\begin{equation}J( c_g)={\partial ( f_3(c_g(\tilde t_1)),   \cdots,  f_3(c_g(\tilde t_{5d}))\over \partial (\theta_1^0, \cdots, 
\theta_d^2, \gamma_1, \cdots, \gamma_{2d})}
\end{equation}

is equal to 
a diagonal matrix $D$ whose diagonal entries are
\begin{equation}
{\partial f_3(c_g(\tilde t_1))\over \partial \theta_1^0}, \cdots, {\partial f_3(c_g(\tilde t_{3d}))\over \partial \theta_d^2}, 
{\partial f_3(c_g(\tilde t_{3d+1}))\over \partial \gamma_1}, \cdots, 
{\partial f_3(c_g(\tilde t_{5d}))\over \partial \gamma_{2d}}
\end{equation}
which are all non-zeros.\par
(b) For $i=1, \cdots, 5d$,  $l=0, \cdots, 4$
$${\partial f_3(c_g(\tilde t_i))\over \partial r_l}=0.$$

\end{lemma}

\bigskip

\begin{proof}  Note $\mathring{ \theta_j^i}, i=0, \cdots, n, j=1, \cdots, d$ are distinct and $\mathring\gamma_k, k=0, \cdots, 2d$ are also distinct.  Thus the coordinates
in Definition 3.3 exist.   Applying $C_M'$ coordinates to $f_3(c(t))$, we have 
\begin{equation} f_3 (c(t))=r_0r_1 r_{2} R \prod_{i=0, j=1, k=1}^{i={2}, j=d, k=2d }(t-\theta_j^i)(t-\gamma_k).\end{equation}
Notice right hand side of (3.12) is in analytic coordinates $C_M'$, and $R$ is a polynomial in variables $r_1, \cdots, r_4$. 
Both parts of Lemma 3.4  follow from the expression (3.12).
We complete the proof. 
\end{proof}

\bigskip

\subsection{Non vanishing of the differential form }
\bigskip

At last we go back to the beginning to construct $\mathbb L$, then $\omega(\mathbb L, \mathbf t)$.
All ingredients in the hidden key will be carefully examined, and remaining terms
in the specialization will be calculated one-by-one. \bigskip

\begin{proof} of Proposition 1.5: 
We compute the non-vanishing $\omega(\mathbb L, \mathbf t)$ in two steps: (1) reduce it to the square matrix $\mathcal A$, (2) break the $\mathcal A$ and calculate two of its block matrices.\medskip

\bigskip

Step 1:  We'll use a generic  $f_0$.  It suffices to prove it for special choices 
of $f_1, f_2$ and distinct $t_1, \cdots, t_{5d+1}$.  
So let $z_0, \cdots, z_4$ be  general homogeneous coordinates of
$\mathbf P^4$. Let $$f_2=z_0z_1z_2z_3z_4.$$  
Let 
$$f_1=z_0z_1z_2 q,$$
where $q$ is a generic quadratic homogeneous polynomial in $z_0, \cdots, z_4$.  
Such a choice of quintics satisfies the pencil condition in Proposition 1.4.
The pencil condition for a neighborhood $\mathbb L$ of $[f_0]$ is equivalent to saying that any non-zero vector $\overrightarrow g$ of 
the plane spanned by $\overrightarrow f_1, \overrightarrow f_2$ is sent to a section
\begin{equation} <\overrightarrow g>\not\in H^0(T_{c_0/X_0}).\end{equation} 
Because the coordinates in $f_1, f_2$ are generic,  
 the pencil condition, i.e. (3.13)  is satisfied. 
Let $$c_g\in I_{\mathbb L}$$
be a generic point in $I_{\mathbb L}$ where $\mathbb L$ is an sufficiently small analytic open set of the plane 
spanned by $f_0, f_1, f_2$.  By the genericity of $f_0$, 
we may assume 
$$c_g=(c_g^0, \cdots, c_g^4)$$ satisfies that   $c_g^i\neq 0$ for all $i$ and equations
$$c_g^i(t)=0, i=0, \cdots, 4$$ have $5d$ distinct zeros $\mathring{\theta_j^i}\in \mathbf P^1$. 
(it is quite important to notice that $c_g$ does not lie in any individual $f_0, f_1, f_2$, but it does lie in a linear combination of them). 
 To calculate the Jacobian matrix $J_{\mathbb L}$, we need to choose auxiliary data: defining equations of $I_{\mathbb L}$ and local coordinates 
of $M$. To have defining equations, we 
 choose $5d+1$ distinct points $t_i$ on $\mathbb C\subset \mathbf P^1$, denoted by
$\mathbf t'=(t_1, \cdots, t_{5d+1})$.    \par
(1) $t_{5d+1}$ is free and $t_1, t_2$ are general satisfying
\begin{equation}
\left|  \begin{array}{cc} f_2(c_g(t_1)) &  f_1(c_g(t_1))\\
f_2(c_g(t_2)) & f_1(c_g(t_2))\end{array}\right|=0,
\end{equation}

 \par
(2) $t_3, \cdots, t_{5d}$ are the $5d-2$ complex numbers
$$\begin{array}{c}
\mathring{ \theta_j^i}, \mathring{\gamma_k}, \quad (i, j)\neq (0, 1), (1, 1)\\
\\
1\leq k\leq 2d,  0\leq i\leq 2, 1\leq j\leq d.
\end{array}$$
satisfying that  $\mathring{\gamma_k}$ are the zeros of 
\begin{equation}
 \delta_1 q(c_g(t))+ \delta_2 z_3z_4|_{c_g(t)}=0,
\end{equation}
with \begin{equation}\begin{array} {c} \delta_1=\left|  \begin{array}{cc} f_0(c_g(t_1)) &  f_2(c_g(t_1))\\
f_0(c_g(t_2)) & f_2(c_g(t_2))\end{array}\right|, \\ \quad \\
\delta_2=\left|  \begin{array}{cc} f_1(c_g(t_1)) &  f_0(c_g(t_1))\\
f_1(c_g(t_2)) & f_0(c_g(t_2))\end{array}\right|.\end{array}\end{equation}
and $\mathring{\theta_j^i}$ are all zeros of
$$c_g^0(t)c_g^1(t)c_g^2(t)=0, $$
but excluding $ \mathring{\theta_1^0}, \mathring{\theta_1^1}$. 
Hence  $t_3, \cdots, t_{5d}$ are just all zeros of 
\begin{equation}
\delta_1 f_1(c_g(t))+\delta_2 f_2(c_g(t))=0.
\end{equation}
but excluding two zeros $\mathring{\theta_1^0}, \mathring{\theta_1^1}$.  
They are distinct because $c_g$ is generic in $I_{\mathbb L}$.
We claim that 
\begin{equation}  (\delta_1,\delta_2)
\end{equation}is generic in $\mathbb C^2$. 
Proof of the claim: The curve $c_g$ lies in  $I_{\mathbb L}$, but does not lie in 
$$P(\Gamma_{span(f_0, f_1)}), \ and\  P(\Gamma_{span(f_0, f_2)}),$$
where $P$ is the projection $M\times S\to M$.   Hence  vectors
$$\{\biggl (f_0(c_g(t)), f_1(c_g(t))\biggr)\}_{t\in \mathbf P^1}$$ span
$\mathbb C^2$.  This implies $$\delta_2=\left|  \begin{array}{cc} f_1(c_g(t_1)) &  f_0(c_g(t_1))\\
f_1(c_g(t_2)) & f_0(c_g(t_2))\end{array}\right|$$
is generic in $\mathbb C$. 
Similarly  $\delta_1$ is generic in $\mathbb C$. By the generacity of the coordinates (that determine $f_1, f_2$), 
$(\delta_1, \delta_2)\in \mathbb C^2$ is generic.  
\bigskip

These $5d+1$ points give specific generators  of the scheme $I_{\mathbb L}$:
$$f(c(t_1)), \cdots, f(c(t_{5d+1})).$$ 

Then we recall the formulation of differential algebra in the introduction as follows.  
As in the introduction, we obtain the $5d-2$ form $\omega(\mathbb L, \mathbf t')$ (
which is a collection of all maximal minors of the Jacobian matrix). 
As in (1.10),  we expand 1-forms $\phi_i, i=3, \cdots, 5d+1$ to obtain that
\begin{equation}\begin{array}{cc}  \phi_i=\left|  \begin{array}{cc} f_1(c(t_1)) &  f_0(c(t_1))\\
f_1(c(t_2)) & f_0(c(t_2))\end{array}\right| \mathbf d f_2(c(t_i))+\left|  \begin{array}{cc} f_0(c(t_1)) &  f_2(c(t_1))\\
f_0(c(t_2)) & f_2(c(t_2))\end{array}\right| \mathbf d f_1(c(t_i))& \\  +\left|  \begin{array}{cc} f_2(c(t_1)) &  f_1(c(t_1))\\
f_2(c(t_2)) & f_1(c(t_2))\end{array}\right|\mathbf d f_0(c(t_i))+ \sum_{l=0, j=1}^{l=2, j=2} h_{lj}^i(c_g) \mathbf d f_l(c(t_j))
\end{array}\end{equation}

By the assumption for $t_1, t_2$, $$\left|  \begin{array}{cc} f_2(c_g(t_1)) &  f_1(c_g(t_1))\\
f_2(c_g(t_2)) & f_1(c_g(t_2))\end{array}\right|=0.$$
We obtain  
 \begin{equation}\begin{array}{c}  \phi_i|_{c_g}=\delta_1 \mathbf d f_1(c(t_i))+\delta_2 \mathbf d f_2(c(t_i))+ \sum_{l=0, j=1}^{l=2, j=2} h_{lj}^i(c_g) \mathbf d f_l(c(t_j))\\\\
=\mathbf d  f_3(c(t_i))+ \sum_{l=0, j=1}^{l=2, j=2} h_{lj}^i(c_g) \mathbf d f_l(c(t_j))
\end{array}\end{equation}
where $$ f_3=\delta_1 f_1+\delta_2 f_2.$$
Notice $(\delta_1, \delta_2)$ is generic. 

\bigskip

The expression says that non-vanishing  of $\omega (\mathbb L, \mathbf t')$ is 
the linear independence of the $5d+5$ differential 1-forms, 
\begin{equation}\begin{array}{c}\mathbf  d f_3(c(t_3)), \cdots, \mathbf d f_3(c(t_{5d+1}))\\
\mathbf d f_0(c(t_1)), \mathbf d f_1(c(t_1)), \mathbf d f_2(c(t_1)),\\
\mathbf d f_0(c(t_2)), \mathbf d f_1(c(t_2)), \mathbf d f_2(c(t_2))
\end{array}\end{equation}
in the cotangent space $(T_{c_g} M)^\ast$, 
  where 
6 functions
$$\begin{array}{c} f_0(c(t_1)), f_1(c(t_1)), f_2(c(t_1)),\\
f_0(c(t_2)), f_1(c(t_2)), f_2(c(t_2))\end{array}$$ come
form the linear combinations 
$$\sum_{l=0, j=1}^{l=2, j=2} h_{lj}^i(c_g) \mathbf d f_l(c(t_j))$$
in (3.20).  This is clearly equivalent to 
the non-degeneracy of Jacobian matrix

 \begin{equation}
\mathcal A(C_M, f_0, f_1, f_2, \mathbf t')={\partial \biggl(f_3(c(t_3)), \cdots, f_3(c(t_{5d+1}), 
f_0(c(t_1)), \cdots, f_2(c(t_2))\biggr)\over \partial C_M}, 
\end{equation}
at $c_g$, 
where $C_M$ is any analytic coordinates' chart of $M$ around $c_g$.

\bigskip
Then it suffices to prove Proposition 1.6.
\bigskip

\it {
 The  $(5d+5)\times (5d+5)$ matrix \begin{equation}
\mathcal A(C_M, f_0, f_1, f_2, \mathbf t')
\end{equation}
is non-degenerate at generic point $c_g$ of  $I_{\mathbb L}$.}

\bigskip

Step 2: 
\begin{proof} of Proposition 1.6:  
With the above choices, $(\delta_1, \delta_2), q$ are all generic. So we can 
 choose  the quasi-polar coordinates $C_M'$ defined in Definition 3.3 to be the local coordinates around $c_g\in I_{\mathbb L}$.   
 We  divide $\mathcal A$ to  block matrices in the following
\begin{equation}
\left(  \begin{array}{cc} \mathcal A_{11} & \mathcal A_{12}\\
\mathcal A_{21}  & \mathcal A_{22} \end{array}\right)
\end{equation}
 where $\mathcal A_{ij}$ are the  Jacobian matrices:

(a)\begin{equation} \mathcal A_{11}={\partial (f_3(c(t_3)), f_3(c(t_4)), \cdots, f_3(c(t_{5d})))\over
\partial (\theta_2^0, \cdots, \hat \theta_1^1, \cdots, \theta_d^2, \gamma_1, \cdots, \gamma_{2d})}, (   \hat{\cdot}    :=omit)
\end{equation}

(b) \begin{equation}  \mathcal A_{12}={\partial (f_3(c(t_3)), f_3(c(t_4)), \cdots, f_3(c(t_{5d})))\over
\partial (\theta_1^0, \theta_1^1, r_0, r_1, r_2, r_3, R)}
\end{equation}

(c) \begin{equation}  \mathcal A_{21}=\end{equation}
$$
{\partial (f_3(c(t_{5d+1})), f_2(c(t_1)), f_2(c(t_2)), f_1(c(t_1)), f_1(c(t_2)), f_0(c(t_1)), f_0(c(t_2)))\over
\partial (\theta_2^0, \cdots, \hat \theta_1^1, \cdots, \theta_d^2, \gamma_1, \cdots, \gamma_{2d})}$$

(d) \begin{equation}  \mathcal A_{22}=\end{equation}
$$
{\partial (f_3(c(t_{5d+1})), f_2(c(t_1)), f_2(c(t_2)), f_1(c(t_1)), f_1(c(t_2)), f_0(c(t_1)), f_0(c(t_2)))\over
\partial (\theta_1^0, \theta_1^1, r_0, r_1, r_2, r_3, R)}.$$
($\mathcal A_{22}$ is a $7\times 7$ matrix.).    

Applying the preparation in sections 3.1, 3.2, we can  found that   $\mathcal A_{11}|_{c_g}$ is a non-zero diagonal matrix and $$\mathcal A_{12}|_{c_g}=0.$$
(Using Lemma 3.4 ). 
Therefore it suffices to show 
\begin{equation}
det(\mathcal A_{22}|_{c_g})\neq 0.
\end{equation}
 Notice $t_{5d+1}$ is generic on $\mathbf P^1$. The genericity of $q$ makes
curve in $\mathbb C^7$, 
\begin{equation}
({\partial f_3(c(t))\over \partial \theta_1^0}, {\partial f_3(c(t))\over \partial \theta_1^1}, {\partial f_3(c(t))\over \partial r_0}, \cdots, {\partial f_3(c(t))\over \partial r_4})
\end{equation}
 span the entire space $\mathbb C^7$. This means the first row vector of $$\mathcal A_{22}|_{c_g}$$
is linearly independent of  other 6 row vectors. 
Hence it suffices for us to show the $6\times 6$ Jacobian matrix 
\begin{equation}
\mathcal B(c_g)=\\
\end{equation}
$${\partial(f_2(c(t_1)), f_2(c(t_2)), f_1(c(t_1)), f_1(c(t_2)), f_0(c(t_1)), f_0(c(t_2)))\over \partial 
(\theta_1^0, \theta_1^1, r_1, r_2, r_3, r_4)}\bigg|_{c_g}.$$
 is non degenerate (the column of partial derivatives with respect to $r_0$ is eliminated).
 To attack  $\mathcal B(c_g)$, we continue to specialize. This time we change the point $c_g$. 
To show $\mathcal B(c_g)$ is non-degenerate,  it suffices to show it is
non-degenerate for a special $c_g'\in I_{\mathbb L}$.  To do that, we let $\mathbb L_2$  be an  open set of pencil through $f_0, f_2$.
Then $I_{\mathbb L}$ must contain a component  $I_{\mathbb L_2}$ where $q$ is generic. 

Let $c_g'$ be a generic point  of $I_{\mathbb L_2}$ ($c_g'$ lies in a lower dimensional subvariety $I_{\mathbb L_2}$, but it
is still in $M_d$ because $f_0$ is generic in $S$).
Because $q$ is generic with respect to  1st, 2nd, 5th and 6th rows and $c_g'(t_1), c_g'(t_2)$ are distinct, two middle rows of
the matrix  $\mathcal B(c_g)$, 
 \begin{equation}\begin{array}{c}
({\partial f_1(c(t_1))\over \partial \theta_1^0}, {\partial f_1(c(t_1))\over \partial \theta_1^1}, {\partial f_1(c(t_1))\over \partial r_1}, \cdots,   {\partial f_1(c(t_1))\over \partial r_4} )|_{c_g'}\\
({\partial f_1(c(t_2))\over \partial \theta_1^0}, {\partial f_1(c(t_2))\over \partial \theta_1^1}, {\partial f_1(c(t_2))\over \partial r_1}, \cdots,   {\partial f_1(c(t_2))\over \partial r_4} )|_{c_g'}
\end{array}\end{equation}
 in $\mathbb C^6$ must be linearly independent of  1st, 2nd, 5th and 6th rows (because  $q$ can  vary freely as $c_g'$ stays fixed). Then we reduce the non-degeneracy of 
$\mathcal B(c_g')$  to the non-degeneracy of $4\times 4$ matrix 
\begin{equation}
Jac(f_0, c_g')={\partial \bigl(f_2(c(t_1)), f_2(c(t_2)), f_0(c(t_1)), f_0(c(t_2))\bigr)\over \partial (\theta_1^0, r_2, r_3, r_4)}|_{c_g'}.
\end{equation}
  Finally we write down the matrix $Jac(f_0, c_g')$,
\begin{equation} \begin{array}{c} Jac(f_0, c_g')\\
\|\\
\lambda \left(\begin{array}{cccccc}
{1\over t_1-\mathring{\theta_1^0}} &1  & 1 & 1  \\
{1\over t_2-\mathring{\theta_1^0}} &1 & 1& 1 \\
 {\partial  f_0(c_g' (t_1))\over \partial \theta_0^1}  & (z_2{\partial f_0 \over \partial z_2})|_{c_g'(t_1)} &
 (z_{3}{\partial f_0 \over \partial z_{3}})|_{c_g'(t_1)}& (z_4{\partial f_0 \over \partial z_4})|_{c_g'(t_1)}\\
{\partial f_0(c_g' (t_2))\over \partial \theta_0^1} & (z_2{\partial f_0 \over \partial z_2})|_{c_g'(t_2)} &
 (z_{3}{\partial f_0 \over \partial z_{3}})|_{c_g'(t_2)}& (z_4{\partial f_0 \over \partial z_4})|_{c_g'(t_2)}
\end{array}\right), \end{array}\end{equation}
where $\lambda $ is a non-zero complex number.
We further compute to have
\begin{equation}\begin{array}{c}  Jac(f_0, c_g')\\
\|\\
\lambda ({1\over t_1-\mathring{\theta_0^1}}-{1\over t_2-\mathring{\theta_0^1}} )\left (\begin{array}{ccc}
1&1&1\\  (z_2{\partial  f_0\over \partial z_2})|_{c_g'(t_1)} &
(z_3{\partial  f_0\over \partial z_3})|_{c_g'(t_1)} &(z_4{\partial  f_0\over \partial z_4})|_{c_g'(t_1)} \\
(z_2{\partial  f_0\over \partial z_2})|_{c_g'(t_2)} &
(z_3{\partial  f_0\over \partial z_3})|_{c_g'(t_2)} &(z_4{\partial  f_0\over \partial z_4})|_{c_g'(t_2)} 
\end{array}\right), \end{array}
\end{equation}
where $\mathring{\theta_0^1} $ is one of complex roots of $c_g'(t)=0$. 
Since $t_1, t_2$ are only required to satisfy one equation (3.14), by the genericity of $q$, 
we may assume $(t_1, t_2)\in \mathbb C^2$ is
generic.   To  show the non-degeneracy of $Jac(f_0, c_g')$,
we consider the Jacobian, i.e. the determinant of the Jacobian matrix, 

$$J(f)=\left|\begin{array}{ccc}
1&1&1\\  (z_2{\partial  f\over \partial z_2})|_{c_g'(t_1)} &
(z_3{\partial  f\over \partial z_3})|_{c_g'(t_1)} &(z_4{\partial  f\over \partial z_4})|_{c_g'(t_1)} \\
(z_2{\partial  f\over \partial z_2})|_{c_g'(t_2)} &
(z_3{\partial  f\over \partial z_3})|_{c_g'(t_2)} &(z_4{\partial  f\over \partial z_4})|_{c_g'(t_2)} 
\end{array}\right|.
$$
where the image $C_g'=c_g'(\mathbf P^1)$ lies in $div(f)$.  We calculate 
$$
J(f)=\left|\begin{array}{cc}
f_4|_{c_g'(t_1)} & f_5|_{c_g'(t_1)} \\
f_4|_{c_g'(t_2)} &
f_5|_{c_g'(t_2)} 
\end{array}\right|.$$
where $$\begin{array} {c} f_4=z_2{\partial  f\over \partial z_2}-z_4{\partial  f\over \partial z_4}\\
f_5=z_3{\partial  f\over \partial z_3}-z_4{\partial  f\over \partial z_4}
\end{array}$$
are two quintics.  Hence  $J(f)=0$ implies
$$({f_4,  f_5})\bigg|_{c_g'(t_1)}$$ $${(f_4, f_5})\bigg|_{c_g'(t_2)}$$ 
are linearly dependent.  Since $(t_1, t_2)\in \mathbb C^2$ is generic, 
then there exist two complex numbers $\epsilon_1, \epsilon_2$ that are not all zero such that 
$$(c_g')^\ast (\epsilon_1 f_4+\epsilon_2 f_5)=0. $$  
Hence $$\eta=(0, 0, \epsilon_1 c_g'^{2}, \epsilon_2 c_g'^{3}, -\epsilon_1 c_g'^{4}-\epsilon_2 c_g'^{4})$$
in $T_{c_g'}M$ gives a holomorphic section of the normal sheaf $N_{div(f)/c_g'}$, where $c_g'^{i}$ is the $i$-th component of
$c_g'$.  Notice two vanishing components of $\eta$ can be moved to a general position (because
the coordinates of $\mathbf P^4$ is generic), while the normal sheaf $N_{div(f)/c_g'}$ containing $\eta$ 
is fixed (independent of coordinates of $\mathbf P^4$).  This contradiction implies 
$J(f)\neq 0$. Hence $J(f_0)\neq 0$.   
Then the determinant of $Jac(f_0, c_g')$ can't be zero.
We complete the proof of Proposition 1.6, therefore of Proposition 1.5. 

\end{proof}

\end{proof}\bigskip

As indicated in the introduction, Proposition 1.5 implies part (2) of Proposition 1.4. Then the section 2 showed this is 
sufficient for theorem 1.1. 
 
\Addresses
\end{document}